\documentclass[french,openright,12pt,twoside]{amsart}
\usepackage{a4}
\usepackage[T1]{fontenc}
\usepackage[all]{xy}
\usepackage{babel,indentfirst}
\usepackage{verbatim}
\usepackage{stmaryrd}
\usepackage{amsfonts}
\usepackage{amsmath,amscd,amssymb}
\usepackage{amsthm}
\usepackage{enumerate}
\usepackage{xspace}
\usepackage{ulem}
\usepackage{euscript}
\usepackage{epsfig}

\newtheorem{exem}{Exemple}[section]
\newtheorem{theo}{Th{\'e}or{\`e}me}[section]

\newtheorem{rem}{Remarque}[section]
\newtheorem{rems}{Remarques}[section]
\newtheorem{coro}{Corollaire}[section]

\newtheorem{pro}{Proposition}[section]
\input cyracc.def
\DeclareFontFamily{U}{russian}{}
\DeclareFontShape{U}{russian}{m}{n}
{ <5><6> wncyr5
<7><8><9> wncyr7
<10><10.95><12><14.4><17.28><20.74><24.88> wncyr10 }{}
\DeclareSymbolFont{Russian}{U}{russian}{m}{n}
\DeclareSymbolFontAlphabet{\mathcyr}{Russian}
\makeatletter
\let\@math@cyr\mathcyr
\begin{document}
\title[Vari{\'e}t{\'e}s de descente, gerbes et obstruction de Brauer-Manin]{Vari{\'e}t{\'e}s de descente, gerbes et obstruction de Brauer-Manin}
\author[Jean-Claude Douai - Michel Emsalem - St{\'e}phane Zahnd]{Jean-Claude Douai - Michel Emsalem - St{\'e}phane Zahnd $^{\left(\ast\right)}$}\thanks{$^{\left(\ast\right)}$ Universit{\'e} des Sciences et Techniques de Lille I, Laboratoire AGAT, B{\^a}timent M2, 
59655 Villeneuve d'Ascq, France.}
\date{18 mars 2003}
\maketitle
\begin{abstract} Nous montrons comment associer {\`a} une gerbe d{\'e}finie sur un corps de nombres une obstruction de Brauer-Manin mesurant, comme dans le cas des vari{\'e}t{\'e}s, le d{\'e}faut d'existence d'une section globale. Ceci nous conduit {\`a} une g{\'e}n{\'e}ralisation de la dualit{\'e} de Tate-Poitou au cas non-ab{\'e}lien.
\end{abstract}
\thispagestyle{empty}
\vspace{5mm}
\tableofcontents
\vspace{1mm}
\section*{Introduction}
Soient $k$ un corps de nombres, $\bar{f}$ un $\bar{k}$-rev{\^e}tement de corps des modules $k$, et $\mathcal{G}\left(\bar{f}\right)$ la $k$-gerbe associ{\'e}e {\`a} $\bar{f}$ (\textit{i.e} la gerbe des mod{\`e}les de $\bar{f}$, cf \cite{DD1}). Dans \cite{DDM}, nous avons introduit la notion de vari{\'e}t{\'e} de descente associ{\'e}e {\`a} $\bar{f}$. Si $V$ est une telle vari{\'e}t{\'e}, la gerbe $\mathcal{G}\left(\bar{f}\right)$ est alors isomorphe au champ quotient $\left[V/GL_{n}\right]$, pour un $n$ idoine. En fait, il existe une infinit{\'e} possible de telles vari{\'e}t{\'e}s de descente $V$ correspondant {\`a} une infinit{\'e} de choix possibles pour l'entier $n$. Soit maintenant $K$ une extension de $k$: tout $K$-point de $V$ d{\'e}finit un $K$-point\footnote{Par $K$-point de $\mathcal{G}$, on entend section de la gerbe au-dessus de $Spec\:K$, ou encore objet de la cat{\'e}gorie fibre de $\mathcal{G}$ au-dessus de $Spec\:K$ (que nous noterons $\mathcal{G}\left(Spec\:k\right)$).} de $\mathcal{G}$, et r{\'e}ciproquement tout $K$-point de $\mathcal{G}$ se rel{\`e}ve en un $K$-point de $V$. Il s'ensuit que si $V$ et $V'$ sont deux vari{\'e}t{\'e}s de descente correspondant {\`a} la m{\^e}me $k$-gerbe $\mathcal{G}$, alors:
	\[V\left(K\right)\neq\emptyset\Leftrightarrow V'\left(K\right)\neq\emptyset
\]
Forts de ces observations, on veut comparer les invariants de $V$ et $V'$; ils ne d{\'e}pendent que de $\mathcal{G}$. En particulier:
	\[Br_{a}V\approx Br_{a}V'\approx Br_{a}\mathcal{G},
\]
et
	\[Pic\:V\approx Pic\:V'\approx Pic\:\mathcal{G}.
\]
Ceci nous am{\`e}ne {\`a} calculer l'invariant de Brauer-Manin\footnote{Le
\textquotedblleft{$\EuScript{H}$}\textquotedblright\ de
\textquotedblleft{$m_{\EuScript{H}}\left(V\right)$}\textquotedblright\ est mis pour faire r{\'e}f{\'e}rence au Principe de Hasse.} $m_{\EuScript{H}}\left(V\right)$ de $V$, {\`a} introduire l'invariant de Brauer-Manin $m_{\EuScript{H}}\left(\mathcal{G}\right)$ de la gerbe $\mathcal{G}$, puis {\`a} prouver que $m_{\EuScript{H}}\left(V\right)=m_{\EuScript{H}}\left(\mathcal{G}\right)$, l'int{\'e}r{\^e}t de cette {\'e}galit{\'e} {\'e}tant sa validit{\'e} pour toute vari{\'e}t{\'e} de descente $V$ correspondant {\`a} $\mathcal{G}$ (plus loin, nous dirons que $V$ est une \textit{pr{\'e}sentation} de $\mathcal{G}$).

Tout ce qui pr{\'e}c{\`e}de s'{\'e}tend aux $k$-gerbes quelconques localement li{\'e}es par un groupe fini (pour des raisons {\'e}videntes, de telles gerbes seront appel{\'e}es \textit{gerbes de Deligne-Mumford}). L'application $m_{\EuScript{H}}$ qui {\`a} une classe de $k$-gerbes $\left[\mathcal{G}\right]$ associe l'invariant de Brauer-Manin d'un de ses repr{\'e}sentants peut alors {\^e}tre vue comme une g{\'e}n{\'e}ralisation de la dualit{\'e} de Tate-Poitou dans le cas ab{\'e}lien (nous renvoyons au th{\'e}or{\`e}me 4.1 pour un {\'e}nonc{\'e} pr{\'e}cis); cet invariant vit dans le groupe de Tate-Shafarevich
	\[\mathcyr{SH}^{1}\left(k,\widehat{\bar{H}}\right)^{D}
\]
o{\`u} $\bar{H}$ est le groupe d'automorphismes d'un objet de $\mathcal{G}\left(Spec\:\bar{k}\right)$.
\vspace{1mm}
\newline

\subsection*{Notations} $k$ d{\'e}signe dans ces lignes un corps de caract{\'e}ristique nulle (souvent un corps de nombres), dont on fixe une cl{\^o}ture alg{\'e}brique $\bar{k}$. Pour toute $k$-vari{\'e}t{\'e} alg{\'e}brique $V$, nous appelons groupe de Brauer cohomologique (ou simplement groupe de Brauer, lorsqu'aucune confusion n'est possible) de $V$ et nous noterons $Br\:V$ le groupe $H^{2}_{\acute{e}t}\left(V,\mathbb{G}_{m}\right)$. Nous noterons $\bar{V}$ la $\bar{k}$-vari{\'e}t{\'e} obtenue {\`a} partir de $V$ par extension des scalaires, le carr{\'e} suivant {\'e}tant alors commutatif:
	\[\xymatrix{\bar{V}=V\otimes_{k}\bar{k} \ar[d] \ar@//[r]^{\epsilon} & V \ar[d]^{p}\\ Spec\:\bar{k} \ar[r]& Spec\:{k}}
\]
De ce diagramme, on d{\'e}duit deux morphismes:
	\[Br\:k\stackrel{\widetilde{p}}{\longrightarrow}Br\:V\stackrel{\widetilde{\epsilon}}{\longrightarrow}Br\:\bar{V}
\]
On note $Br_{0}V$ (resp. $Br_{1}V$) l'image de $\widetilde{p}$ (resp. le noyau de $\widetilde{\epsilon}$), et $Br_{a}V$ le quotient $Br_{1}V/Br_{0}V$. Pour tout groupe alg{\'e}brique $G$, $\widehat{G}$ d{\'e}signe le groupe des caract{\`e}res de $G$; lorsque $A$ est un groupe ab{\'e}lien, nous noterons $A^{D}$ le dual (de Pontrjagin) de $A$: $Hom\left(A,\mathbb{Q}/\mathbb{Z}\right)$. Si $\bar{H}$ est un $Gal\left(\bar{k}/k\right)$-module, on note:
	\[H^{i}\left(k,\bar{H}\right)=H^{i}\left(Gal\left(\bar{k}/k\right),\bar{H}\right)\ \ \ \left(i=1,2\right)
\]
Si de plus $k$ est un corps de nombres, on d{\'e}finit:
	\[\mathcyr{SH}^{i}\left(k,\bar{H}\right)=\ker\left\{H^{i}\left(k,\bar{H}\right)\longrightarrow\prod_{all\:v}{H^{i}\left(k_{v},\bar{H}\right)}\right\}\ \ \ \left(i=1,2\right)
\]
\begin{section}{Rappels}
\begin{subsection}{Calcul de $Br_{a}V$ dans le cas o{\`u} $V$ est un espace homog{\`e}ne de $SL_{n}$ avec isotropie $H$}
$\ $
\newline

Soient $k$ un corps de caract{\'e}ristique nulle, et $V$ une $k$-vari{\'e}t{\'e} alg{\'e}brique lisse, g{\'e}om{\'e}triquement irr{\'e}ductible. De la suite spectrale
	\[H^{p}\left(k,H^{q}_{\acute{e}t}\left(\bar{V},\mathbb{G}_{m}\right)\right)\Longrightarrow H^{p+q}_{\acute{e}t}\left(V,\mathbb{G}_{m}\right)
\]
on d{\'e}duit la suite exacte longue
\begin{equation}
	\xymatrix{0 \ar[r] & H^{1}\left(k,\bar{k}\left[V\right]^{\ast}\right) \ar[r] & Pic\:V \ar[r] & Pic\:\bar{V}^{Gal\left(\bar{k}/k\right)} \ar`dr_l[ll]`^dr[ll] [dl]  & \\ &&Br_{1}V \ar[r]&H^{1}\left(k,Pic\:\bar{V}\right) \ar[r] & H^{3}\left(k,\bar{k}\left[V\right]^{\ast}\right)}
\end{equation}
Posons:
	\[U\left(\bar{V}\right)=\frac{\bar{k}\left[V\right]}{\bar{k}^{\ast}}
\]
La suite exacte $\left(1\right)$ fournit une nouvelle suite exacte:
\begin{equation}
	Pic\:\bar{V}^{Gal\left(\bar{k}/k\right)}\rightarrow H^{2}\left(k,U\left(\bar{V}\right)\right) \rightarrow Br_{a}V \rightarrow H^{1}\left(k,Pic\:\bar{V}\right) \rightarrow H^{3}\left(k,U\left(\bar{V}\right)\right)
\end{equation}
Supposons que $V$ est un $k$-espace homog{\`e}ne d'un $k$-groupe alg{\'e}brique semi-simple simplement connexe $\widetilde{G}$ (\textit{e.g} $SL_{n}$) avec isotropie un groupe fini, c'est-{\`a}-dire: il existe un $\bar{k}$-groupe fini $\bar{H}$ tel que:
	\[\bar{V}=\widetilde{G}\left(\bar{k}\right)/\bar{H}
\]
Nous avons alors la suite exacte 
	\[0\longrightarrow U\left(\bar{V}\right) \longrightarrow U\left(\widetilde{G}\left(\bar{k}\right)\right)
\]
provenant de la fibration $\widetilde{G}\left(\bar{k}\right)\rightarrow \bar{V}$. Or on sait (cf le lemme 6.5 (iii) de \cite{Sa}) que:

\[U\left(\widetilde{G}\left(\bar{k}\right)\right)=\widehat{\widetilde{G}\left(\bar{k}\right)}=0
\]
La suite exacte $\left(2\right)$ se r{\'e}duit alors {\`a} l'isomorphisme:
	\[Br_{a}V\stackrel{\sim}{\longrightarrow}H^{1}\left(k,Pic\:\bar{V}\right)
\]
\begin{rem} \textup{Notons au passage que cet isomorphisme tient encore lorsque $V$ est une vari{\'e}t{\'e} alg{\'e}brique propre (\textit{e.g} projective) d{\'e}finie sur un corps de nombres $k$. Car dans cette situation, d'une part le groupe $H^3\left(k,\mathbb{G}_{m}\right)$ est nul, et d'autre part $\bar{k}\left[V\right]^{\ast}$ se r{\'e}duit {\'e}videmment aux constantes.}

\textup{Par exemple, ($k$ {\'e}tant toujours un corps de nombres) le groupe $Br_{a}V$ (donc \textit{a fortiori} $\mathcyr{B}\left(V\right)$) est nul lorsque $V$ est une $k$-vari{\'e}t{\'e} de Severi-Brauer, ou une $k$-vari{\'e}t{\'e} projective lisse qui est une intersection compl{\`e}te de dimension $\geq3$. En effet dans les deux cas, on a: $Pic\:\bar{V}=\mathbb{Z}$ (c'est {\'e}vident pour les vari{\'e}t{\'e}s de Severi-Brauer, et c'est essentiellement une cons{\'e}quence du th{\'e}or{\`e}me de la section hyperplane de Lefschetz dans le deuxi{\`e}me cas).}
\end{rem}
\end{subsection}
\begin{subsection}{Exemples}
\begin{enumerate}[(i)]
\item Si $H=0$, alors $\bar{V}\approx\widetilde{G}\left(\bar{k}\right)$, et $Pic\:\bar{V}=0$ (car $Pic\:\widetilde{G}\left(\bar{k}\right)=0$ par le corollaire 4.5 de \cite{FI}).
\item Si $H=\mu$ est un $k$-sous-groupe central de $\widetilde{G}$, alors $V=G=\widetilde{G}/\mu$ est semi-simple, et $Pic\:\widetilde{G}\left(\bar{k}\right)=\widehat{\mu\left(\bar{k}\right)}$ (par le corollaire 4.6 de \cite{FI}), d'o{\`u}:

\[Br_{a}V=Br_{a}G=H^{1}\left(\bar{k}/k,\widehat{\mu\left(\bar{k}\right)}\right)=H^{1}\left(k,\widehat{\mu}\right)
\]
$Pic\:G$ et $Br_{a}G$ sont justiciables de la philosophie de Kottwitz: ce sont des invariants des groupes semi-simples qui sont nuls lorsque $G=\widetilde{G}$ est simplement connexe. Ils peuvent donc s'exprimer en fonction du centre $Z\left(^{L}G\right)$ du dual de Langlands de $G$ \cite{K}; lorsque $G$ est semi-simple, ce dernier coïncide avec le dual du noyau du rev{\^e}tement universel de $G$.

Cette remarque vaut encore pour
	\[\mathcyr{B}\left(G\right)=\ker\left\{Br_{a}G\rightarrow{\prod_{all\ v}Br_{a}G_{v}}\right\}
\]
dans le cas o{\`u} $k$ est un corps de nombres (le produit {\'e}tant pris sur toutes les places $v$ de $k$).
\item Soient $k$ un corps de nombres et
prenons pour $V$ un $k$-tore $T$. Alors (cf le lemme 6.9 de \cite{Sa}):
	\[Pic\:\bar{T}=H^{1}\left(k,\widehat{T}\right)\ \ et\ \ Br_{a}T=H^{2}\left(k,\widehat{T}\right)
\]
En outre:
	\[\mathcyr{B}\left(T\right)\approx\mathcyr{SH}^{2}\left(k,\widehat{T}\right)\approx\mathcyr{SH}^{1}\left(k,\widehat{T}\right)^{D}
\]
le deuxi{\`e}me isomorphisme {\'e}tant fourni par la dualit{\'e} de Kottwitz \cite{K} qui {\'e}tend aux tores celle de Tate-Poitou.
\item Consid{\'e}rons maintenant un $k$-espace homog{\`e}ne de $SL_{n}$ avec isotropie un groupe fini; on suppose donc qu'il existe un groupe fini $\bar{H}$ tel que:
	\[\bar{V}\approx SL_{n,\bar{k}}/\bar{H}
\]
Alors $Pic\:\bar{V}\approx\widehat{\bar{H}}$ (cf \cite{BK2}). On dispose en effet de la $\bar{k}$-fibration:
	\[SL_{n,\bar{k}}\longrightarrow SL_{n,\bar{k}}/\bar{H}
\]
{\`a} laquelle est attach{\'e}e la suite spectrale
	\[E^{p,q}_{2}=H^{p}\left(\bar{H},H^{q}\left(SL_{n},\mathbb{G}_{m}\right)\right)\Longrightarrow H^{p+q}_{\acute{e}t}\left(\bar{V},\mathbb{G}_{m}\right)
\]
Dans cette derni{\`e}re, le terme $E_{2}^{0,1}$ est nul\footnote{En effet,
$Pic\:SL_{n}=H^{1}\left(SL_{n,\bar{k}},\mathbb{G}_{m}\right)=Hom\left(\Pi_{1}\left(SL_{n}\right),\mathbb{G}_{m}\right)=0$ puisque $SL_{n}$ est simplement connexe.}, donc:
\[H^{1}_{\acute{e}t}\left(\bar{V},\mathbb{G}_{m}\right)=H^{1}\left(\bar{H},\mathbb{G}_{m}\right)=Hom\left(\bar{H},\mathbb{G}_{m}\right)=\widehat{\bar{H}}
\]
D'o{\`u} la:
\end{enumerate}
\begin{pro} Soit $V$ un $k$-espace homog{\`e}ne d'un groupe semi-simple
simplement connexe $\widetilde{G}$ avec isotropie un groupe fini $\bar{H}$.
Alors \cite{BK2}:
	\[Br_{a}V\approx H^{1}\left(k,\widehat{\bar{H}}\right)
\]
En outre, si $k$ est un corps de nombres, et si on suppose que $V$ a des points localement partout (i.e si $V_{v}=V\otimes_{k}k_{v}$ a un $k_{v}$-point, pour toute place $v$ de $k$), alors:
	\[\mathcyr{B}\left(V\right)=\ker\left\{Br_{a}V\rightarrow{\prod_{all\ v}Br_{a}V_{v}}\right\}\approx \mathcyr{SH}^{1}\left(k,\widehat{\bar{H}}\right)
\]
\end{pro}
\begin{coro} Sous les hypoth{\`e}ses et notations de la proposition pr{\'e}c{\'e}dente, si $\bar{H}$ est sans caract{\`e}re, alors $Br_{a}V=\mathcyr{B}\left(V\right)=0$.
\end{coro}
\begin{exem} \textup{Le corollaire s'applique donc si $\bar{H}=SL\left(2,\mathbb{F}_{p}\right)$ avec $p\neq2,3$, ou encore si $\bar{H}=\EuScript{A}_{n}$ avec $n\geq5$. Plus g{\'e}n{\'e}ralement, il suffit que $\bar{H}$ soit {\'e}gal {\`a} son groupe d{\'e}riv{\'e}.}
\end{exem}
\end{subsection}
\end{section}
\begin{section}{Interpr{\'e}tation comme champ quotient des $k$-gerbes localement li{\'e}es par un groupe alg{\'e}brique fini}
On s'int{\'e}resse donc dans cette section aux $k$-gerbes qui sont des champs de Deligne-Mumford \cite{LMB}. Rappelons d'abord la proposition 5.1 de \cite{DDM}:
\begin{pro} Soient $k$ un corps, et $\mathcal{G}$ une $k$-gerbe (pour la topologie {\'e}tale) qui est un champ
de Deligne-Mumford. Alors:
\begin{enumerate}
\item Il existe une $k$-alg{\`e}bre $L$ avec action {\`a} gauche d'un groupe fini $\Gamma$ admettant $k$ comme anneau des invariants telle que $\mathcal{G}$ soit isomorphe au champ quotient $\left[Spec\:L/\Gamma\right]$;
\item Il existe un $k$-sch{\'e}ma affine $V$, un entier $n\geq0$, une action {\`a} droite de $GL_{n,k}$ sur $V$ et un $1$-morphisme $\pi:V\rightarrow\mathcal{G}$ avec les propri{\'e}t{\'e}s suivantes:
\begin{enumerate}[(i)]
\item $\pi$ induit un isomorphisme du champ quotient $\left[V\right/GL_{n,k}]$ vers $\mathcal{G}$;
\item $V$ est lisse et g{\'e}om{\'e}triquement irr{\'e}ductible;
\item l'action de $GL_{n,k}$ sur $V$ est transitive et {\`a} stabilisateurs finis;
\item pour chaque extension $K$ de $k$, chaque objet de $\mathcal{G}\left(K\right)$ se rel{\`e}ve en un point de $V\left(K\right)$ via $\pi$.
\end{enumerate}
En particulier, {\`a} cause de (iii) et (iv), si $K$ est une extension de $k$ telle que
$\mathcal{G}\left(K\right)\neq\emptyset$,\footnote{Ce qui signifie que l'ensemble d'objets
$Ob\left(\mathcal{G}\left(K\right)\right)$ de la cat{\'e}gorie fibre de $\mathcal{G}$ au-dessus de $Spec\:K$
est non-vide; par la suite, nous fairons syst{\'e}matiquement cet abus de langage.} la $K$-vari{\'e}t{\'e}
$V\otimes_{k}K$ est isomorphe au quotient de $GL_{n,K}$ par un groupe fini.
\end{enumerate}
\end{pro}
\begin{rems}$\ $

\begin{enumerate}[(a)]
\item \textup{Dans la remarque 5.2(b) de \cite{DDM}, il est montr{\'e} que l'on peut en fait prendre pour $k$-alg{\`e}bre $L$ une extension galoisienne finie de $k$, auquel cas $\Gamma$ est l'ensemble des couples $\left(\sigma,\varphi\right)$ o{\`u} $\sigma\in\:Gal\left(L/k\right)$ et $\varphi:\sigma x\stackrel{\sim}{\rightarrow}x$ est un isomorphisme dans la cat{\'e}gorie (en fait, le groupoïde) $\mathcal{G}\left(L\right)$. Il y a une structure de groupe sur $\Gamma$ pour laquelle la projection naturelle $\Gamma\rightarrow Gal\left(L/k\right)$ est un morphisme surjectif, et le noyau $H=H\left(L\right)$ est le stabilisateur fini dont l'existence est donn{\'e}e par le (iii) de la proposition 2.1.}

\textup{A la gerbe $\mathcal{G}$ est aussi associ{\'e}e une extension (cf \cite{Gi}, \cite{Sp}) 
	\[\mathcal{G}\leftrightsquigarrow\left(\EuScript{E}\right):1\rightarrow{H}\rightarrow\Gamma\rightarrow{Gal\left(L/k\right)}\rightarrow{1}
\]
d{\'e}finissant une action ext{\'e}rieure $\mathcal{L}_{H}$ de $Gal\left(L/k\right)$ sur $H=H\left(L\right)$, et une classe de 2-cohomologie not{\'e}e $\left[\mathcal{G}\right]$ dans $H^{2}\left(L/k,\mathcal{L}_{H}\right)\hookrightarrow H^{2}\left(k,\mathcal{L}_{H}\right)$.}
\item \textup{On obtient les m{\^e}mes conclusions en rempla\c{c}ant dans la proposition 2.1 $GL_{n}$ par $SL_{n}$ (cf remarque 5.2(c) de \cite{DDM}).}
\end{enumerate}
\end{rems}
$\ $
\vspace{1mm}
\newline
\begin{center}
\uline{Une construction fondamentale}
\end{center}
$\ $
\newline
Partons de l'extension $\left(\EuScript{E}\right)$ de la remarque pr{\'e}c{\'e}dente:
	\[\left(\EuScript{E}\right):1\rightarrow{H}\rightarrow\Gamma\rightarrow{Gal\left(L/k\right)}\rightarrow{1}
\]
$\Gamma$ est un groupe fini; on peut donc le plonger dans $SL_{n}$ pour un certain $n$, ce qui conduit au diagramme suivant:
	\[\left(D\right):\xymatrix{1\ar[r] & H \ar[r] \ar@{=}[d] & \Gamma \ar[r] \ar[d] & Gal\left(L/k\right) \ar[r] \ar@{-->}[d] & 1\\ 1 \ar[r] & H \ar[r] & SL_{n,\bar{k}} \ar@{-->}[r] & SL_{n,\bar{k}}/H \ar[r] & 1}
\]
$SL_{n,\bar{k}}/H$ n'est pas un groupe, puisque $H$ n'est pas n{\'e}cessairement normal dans $SL_{n,\bar{k}}$. C'est seulement un $k$-espace homog{\`e}ne (toujours au sens de Springer \cite{Sp}), d'o{\`u} la pr{\'e}sence des pointill{\'e}s dans le diagrammme pr{\'e}c{\'e}dent. La fl{\`e}che verticale
	\[\xymatrix{Gal\left(L/k\right) \ar@{-->}[d]\\ SL_{n,\bar{k}}/H }
\]
donne lieu {\`a} un $1$-cocycle dans $Z^{1}\left(L/k;SL_{n},H\right)$, qui repr{\'e}sente pr{\'e}cis{\'e}ment la classe du $k$-espace homog{\`e}ne $V$ du (2) de la proposition 2.1. La $k$-gerbe $\mathcal{G}\approx\left[V/SL_{n}\right]$ (associ{\'e}e {\`a} $\left(\EuScript{E}\right)$) s'interpr{\`e}te alors comme la gerbe des rel{\`e}vements du $k$-espace homog{\`e}ne $V$ {\`a} $SL_{n}$. En d'autres termes, $\left[\mathcal{G}\right]$ est l'image de $V$ par le cobord (cf \cite{Sp}, \cite{D1})
	\[\delta^{1}:Z^{1}\left(L/k;SL_{n},H\right)\longrightarrow H^{2}\left(k,\mathcal{L}_{H}\right)
\]
Dans la suite, nous appellerons \textbf{pr{\'e}sentation de $\mathcal{G}=\left[V/SL_{n}\right]$} un couple $\left(V,\pi\right)$ comme dans la proposition 2.1.
\setcounter{rem}{1}
\begin{rem} \textup{La proposition 2.1 traduit en particulier le fait que les deux groupes de Brauer d'un sch{\'e}ma (cohomologique et \textquotedblleft{Azumaya}\textquotedblright) coïncident lorsque ce sch{\'e}ma est un corps. En effet, si $\mathcal{G}\in H^{2}\left(k,\mathbb{G}_{m}\right)$, alors il existe un $n$ tel que $\mathcal{G}\in H^{2}\left(k,\mu_{n}\right)$, puisque le groupe de Brauer d'un corps, et plus g{\'e}n{\'e}ralement d'un sch{\'e}ma r{\'e}gulier \cite{G2}, est de torsion. Par cons{\'e}quent, $\mathcal{G}$ est un champ de Deligne-Mumford, et il existe un entier $n'$ et un $k$-espace homog{\`e}ne $V$ de $SL_{n'}$ avec isotropie $\mu_{n}$ tels que 
	\[\mathcal{G}\approx\left[V/SL_{n'}\right]
\]
Remarquons que l'entier $n'$ peut {\^e}tre diff{\'e}rent de $n$, et c'est en particulier le cas pour le corps $k_{M}$ construit par Merkuriev dans son article sur la conjecture de Kaplansky \cite{Me}: il existe un {\'e}l{\'e}ment de $H^{2}\left(k_{M},SL_{2}\right)$ qui ne provient pas d'un {\'e}l{\'e}ment de $H^{1}\left(k_{M},PGL_{2}\right)=H^{1}\left(k_{M};SL_{2},\mu_{2}\right)$ (mais qui est atteint par un {\'e}l{\'e}ment de $H^{1}\left(k_{M},PGL_{4}\right)$).}

\textup{Cependant, l'espace homog{\`e}ne $V$ n'est autre que la vari{\'e}t{\'e} de Severi-Brauer pr{\'e}-image de $\mathcal{G}$ par le morphisme naturel
	\[Br_{Az}k\longrightarrow H^{2}\left(k,\mathbb{G}_{m}\right)
\]
et on retrouve ainsi le point de vue de \cite{EHKV}.}
\end{rem}
\end{section}
\begin{section}{Invariant de Brauer-Manin d'une $k$-gerbe localement li{\'e}e par un groupe fini}
Les $k$-champs alg{\'e}briques (en particulier les $k$-gerbes qui sont des champs de Deligne-Mumford) sont des g{\'e}n{\'e}ralisations de la notion de sch{\'e}ma\footnote{Plus rigoureusement, il existe un foncteur pleinement fid{\`e}le de la cat{\'e}gorie des $k$-sch{\'e}mas dans celle des $k$-champs alg{\'e}briques, d{\'e}fini en associant {\`a} un $k$-sch{\'e}ma $S$ le champ discret $S^{ch}$ qu'il repr{\'e}sente. Les objets de $S^{ch}$ n'ont pas d'automorphisme non-trivial, et r{\'e}ciproquement tout $k$-champ alg{\'e}brique dont les objets n'ont pas d'automorphisme non-trivial provient d'un $k$-espace alg{\'e}brique via ce foncteur. Nous renvoyons {\`a} \cite{LMB} pour plus de d{\'e}tails.}. Par suite, il est tout-{\`a}-fait naturel de d{\'e}finir l'obstruction de Brauer-Manin d'une $k$-gerbe de mani{\`e}re analogue {\`a} celle d'un $k$-sch{\'e}ma.

Soit donc $\mathcal{G}$ une $k$-gerbe, qui est un champ de Deligne-Mumford; on suppose son $k$-lien localement
repr{\'e}sentable par un groupe fini $H$. Le site {\'e}tale de $\mathcal{G}$
est d{\'e}fini au chapitre 12 de \cite{LMB}. Le groupe de Brauer cohomologique
$Br\:\mathcal{G}=H^{2}_{\acute{e}t}\left(\mathcal{G},\mathbb{G}_{m}\right)$ est d{\'e}fini dans \cite{Ve}
(o{\`u}, d'une mani{\`e}re plus g{\'e}n{\'e}rale, la cohomologie d'un topos localement annel{\'e} est
d{\'e}finie). Il existe une suite spectrale (cf \cite{Ve}, prop. 5.3):
	\[H^{p}\left(k,H^{q}_{\acute{e}t}\left(\bar{\mathcal{G}},\mathbb{G}_{m}\right)\right)\Longrightarrow H^{p+q}_{\acute{e}t}\left(\mathcal{G},\mathbb{G}_{m}\right)
\]
qui permet de d{\'e}finir $Br_{0}\mathcal{G}$, $Br_{1}\mathcal{G}$ et $Br_{a}\mathcal{G}$ de mani{\`e}re analogue {\`a} celle de la section 1.

Soit $\left(V,\pi\right)$ une pr{\'e}sentation de $\mathcal{G}$ (d'apr{\`e}s nos conventions, on a donc $\mathcal{G}\approx\left[V/SL_{n}\right]$). Posons:
	\[U\left(\bar{\mathcal{G}}\right)=\frac{\bar{k}\left[\mathcal{G}\right]^{\ast}}{\bar{k}^{\ast}}
\]
On a:
	\[U\left(\bar{\mathcal{G}}\right)\subset U\left(\bar{V}\right) \subset U\left(SL_{n,\bar{k}}\right)=\widehat{SL_{n,\bar{k}}}=0
\]
L'analogue de la suite exacte $\left(2\right)$ associ{\'e}e {\`a} la suite spectrale pr{\'e}c{\'e}dente implique alors:
	\[\begin{array}{rl}Br_{a}\mathcal{G}&\approx H^{1}\left(k,Pic\:\bar{\mathcal{G}}\right)\\ &\approx H^{1}\left(k,Hom\left(\Pi_{1}\left(\bar{\mathcal{G}}\right),\mathbb{G}_{m}\right)\right)\\&\approx H^{1}\left(k,Hom\left(\bar{H},\mathbb{G}_{m}\right)\right)\\\end{array}
\]
car il est bien connu que $\Pi_{1}\left(\bar{\mathcal{G}}\right)=\bar{H}$ (cf \cite{No}), et finalement:
\begin{equation}
	Br_{a}\mathcal{G}\approx H^{1}\left(k,\widehat{\bar{H}}\right)
\end{equation}

Si $K$ est un corps qui est une extension quelconque de $k$, nous pouvons d{\'e}finir un accouplement:
	\[\begin{array}{ccl}\mathcal{G}\left(K\right)\times Br\:\mathcal{G}& \rightarrow & Br\:K\\ \left(x,b\right) & \longmapsto & b\left(x\right)\end{array}
\]
o{\`u} l'image $b\left(x\right)$ peut {\^e}tre interpr{\'e}t{\'e}e de diff{\'e}rentes fa\c{c}ons ($\mathcal{B}$ d{\'e}signe ci-dessous un repr{\'e}sentant de $b$):
\begin{enumerate}[(i)]
\item $b\left(x\right)$ est la gerbe r{\'e}siduelle de $\mathcal{B}$ au point $x$ du champ alg{\'e}brique $\mathcal{G}=\left[V/SL_{n}\right]$;
\item $x$ peut {\^e}tre vu comme une section au dessus de $Spec\:K$ du (1-)morphisme structural $\mathcal{G}\rightarrow{Spec\:k}$; autrement dit, c'est un (1-)morphisme rendant commutatif le diagramme (de morphismes de champs) suivant:
	\[\xymatrix{&\mathcal{G} \ar[d]\\Spec\:K \ar[r] \ar[ur]^{x}&Spec\:k}
\]
$\mathcal{B}$ {\'e}tant un champ (puisque c'est une gerbe) sur $\mathcal{G}$, on peut consid{\'e}rer le champ image inverse $x^{\ast}\mathcal{B}$ de $\mathcal{B}$ par le morphisme $x$; la gerbe $x^{\ast}\mathcal{B}$ ainsi obtenue\footnote{L'image inverse d'une gerbe par un morphisme de champs est toujours une gerbe \cite{Gi}.} correspond exactement {\`a} $b\left(x\right)$; elle est donc obtenue par pull-back {\`a} partir de $x$ et de $\mathcal{B}$:
	\[\xymatrix{&\mathcal{B} \ar@{-}[d]\\x^{\ast}\mathcal{B} \ar@{-}[d] \ar@{--}[ur]^{x^{\ast}}&\mathcal{G} \ar[d]\\Spec\:K \ar[r] \ar[ur]^{x}&Spec\:k}
\]
\item on peut interpr{\'e}ter $\mathcal{B}$ comme une alg{\`e}bre d'Azumaya sur $\mathcal{G}$ (\textit{i.e} un torseur sous $PGL_{n}\left(\mathcal{O}_{\mathcal{G}}\right)$ pour la topologie {\'e}tale sur $\mathcal{G}$); ainsi, la fibre de $\mathcal{B}$ en $x$ est simplement une alg{\`e}bre simple centrale sur $K$, et c'est pr{\'e}cis{\'e}ment l'{\'e}l{\'e}ment $b\left(x\right)$ de $Br\:K$ recherch{\'e}. Notons que dans le cas g{\'e}n{\'e}ral, on ne peut utiliser cette interpr{\'e}tation puisque les deux groupes de Brauer (cohomologique et \textquotedblleft{Azumaya}\textquotedblright) ne coïncident pas n{\'e}cessairement; on sait cependant qu'ils sont {\'e}gaux dans de nombreuses situations (cf \cite{G2}, \cite{Ga}, \cite{Sc}\ldots) et que cette {\'e}galit{\'e} tient en particulier dans le cas des sch{\'e}mas affines \cite{Ga}. Or le groupe de Brauer de $\mathcal{G}$ est reli{\'e} au groupe de Brauer d'une quelconque de ses pr{\'e}sentations $V$; comme un tel $V$ est un espace homog{\`e}ne de $SL_{n}$, $V$ est en particulier affine, ce qui rend l{\'e}gitime notre interpr{\'e}tation.
\end{enumerate}

De plus, toute $K$-section de la gerbe $\mathcal{G}$ (\textit{i.e} tout objet de $\mathcal{G}\left(K\right)$) est un $H$-torseur sur $Spec\:K$ ($\mathcal{G}$ est par d{\'e}finition localement {\'e}quivalente {\`a} la gerbe $Tors\:H$; l'existence d'une $K$-section implique que $\mathcal{G}_{\left|K\right.}$ est {\'e}quivalente {\`a} la gerbe des $H$-torseurs sur $K$). Si on suppose que $K$ est un corps local, on obtient alors l'{\'e}nonc{\'e} suivant:
\begin{pro}[Cas local] Soient $K$ un corps local, $\mathcal{G}$ une $K$-gerbe li{\'e}e par un groupe
ab{\'e}lien fini $H$, et $\left(V,\pi\right)$ une pr{\'e}sentation de $\mathcal{G}$. Le diagramme suivant est commutatif:
	\[\xymatrix@C=2pt@R=25pt{\left(Acc.1\right)&V\left(K\right)\ar@{->>}[d]&\times&Br_{a}V\ar[dr] \\\left(Acc.2\right)&\mathcal{G}\left(K\right)\ar@{->>}[d]&\times&Br_{a}\mathcal{G}\ar @<1ex> [u]^{\approx}\ar[r]&\mathbb{Q}/\mathbb{Z} \\\left(Acc.3\right)& H^{1}\left(K,\bar{H}\right)&\times&H^{1}\left(K,\widehat{\bar{H}}\right) \ar @<1ex> [u]^{\approx} \ar[ur]}
\]
o{\`u}:
\begin{itemize}
\item l'accouplement $\left(Acc.1\right):V\left(K\right)\times Br_{a}V\longrightarrow\mathbb{Q}/\mathbb{Z}$ est d{\'e}fini comme suit: {\`a} un point $x$ de $V\left(K\right)$, et {\`a} une classe $b$ dans $Br_{a}V$, on associe:
	\[\left\langle x,b\right\rangle=\left[s_{x}\left(b\right)\right]_{x}
\]
$s_{x}$ d{\'e}signant la section\footnote{En effet, l'existence d'un point $K$-rationnel entra{\^i}ne que la
suite:\[0\longrightarrow Br\:K\longrightarrow Br_{1}V\longrightarrow Br_{a}V\longrightarrow 0
\]
est scind{\'e}e.} induite par $x$ de la projection canonique $p$ (cf \cite{BK1}):
	\[\xymatrix{Br_{1}V \ar[r]_{p}&Br_{a}V \ar@/_12pt/[l]_{s_{x}}}
\]
\item l'accouplement $\left(Acc.2\right):\mathcal{G}\left(K\right)\times Br_{a}\mathcal{G}\longrightarrow\mathbb{Q}/\mathbb{Z}$ est d{\'e}fini de la m{\^e}me mani{\`e}re que $\left(Acc.1\right)$;
\item l'accouplement $\left(Acc.3\right):H^{1}\left(K,\bar{H}\right)\times H^{1}\left(K,\widehat{\bar{H}}\right)\longrightarrow\mathbb{Q}/\mathbb{Z}$ est l'accouplement de Tate pour les corps locaux.
\end{itemize}
\end{pro}
$\ $
\newline

Lorsque $k$ est un corps de nombres, ce que nous supposons {\`a} partir de maintenant, nous pouvons d{\'e}finir pour toute place $v$ de $k$ l'accouplement:
	\[\begin{array}{ccl} \mathcal{G}\left(k_{v}\right)\times Br_{1}\:\mathcal{G}& \rightarrow & \mathbb{Q}/\mathbb{Z}\\ \left(x,b\right) & \longmapsto & inv_{v}\left(b\left(x\right)\right)\end{array}
\]
o{\`u} comme d'habitude $inv_{v}$ est l'invariant donn{\'e} par la th{\'e}orie du corps de classes, et $b\left(x\right)$ est la classe dans $Br\:k_{v}$ de $\mathcal{B}_{x}$, o{\`u} $\mathcal{B}$ est un repr{\'e}sentant de $b$. Supposons maintenant que $\mathcal{G}\left(k_{v}\right)$ soit non-vide pour toute place $v$ de $k$, et restreignons nous au sous-groupe $\mathcyr{B}\left(\mathcal{G}\right)$ de $Br_{a}\mathcal{G}$ d{\'e}fini par:
	\[\mathcyr{B}\left(\mathcal{G}\right)=\ker\left\{Br_{a}\mathcal{G}\longrightarrow \prod_{all\:v}{Br_{a}\:\mathcal{G}_{\left|k_{v}\right.}}\right\}
\]
Nous d{\'e}finissons ainsi un accouplement:
	\[\begin{array}{rccl}\left\langle.\:,\:.\right\rangle: &
	$\(\displaystyle\prod_{all\:v}{\mathcal{G}\left(k_{v}\right)}\)$ \times
	\mathcyr{B}\left(\mathcal{G}\right)& \longrightarrow & \mathbb{Q}/\mathbb{Z}\\ &
	\left(\left(x_{v}\right)_{v},b\right) & \longmapsto &
	$\(\displaystyle\sum_{all\:v}{inv_{v}\left(\widetilde{b}\left(x_{v}\right)\right)}\)$\end{array}
\]
o{\`u} $\widetilde{b}$ est un relev{\'e} de $b$ dans $Br_{1}\mathcal{G}$. Par analogie avec la d{\'e}finition
usuelle de cet accouplement (cf \cite{Bo1}), $\left\langle
\left(x_{v}\right)_{v},b\right\rangle$ ne d{\'e}pend pas de $\left(x_{v}\right)_{v}$, et $\left\langle
\left(x_{v}\right)_{v},b\right\rangle\neq0$ est une obstruction {\`a} l'existence d'une section $k$-rationnelle
de $Spec\:k$ (\textit{i.e} d'un objet de la cat{\'e}gorie fibre $\mathcal{G}\left(k\right)$). Nous obtenons de cette fa\c{c}on un {\'e}l{\'e}ment bien d{\'e}fini:
	\[m_{\mathcal{H}}\left(\mathcal{G}\right)\in \mathcyr{B}\left(\mathcal{G}\right)^{D}=Hom\left(\mathcyr{B}\left(\mathcal{G}\right),\mathbb{Q}/\mathbb{Z}\right)=Hom\left(\mathcyr{SH}^{1}\left(k,\widehat{\bar{H}}\right),\mathbb{Q}/\mathbb{Z}\right)
\]
puisque: $\mathcyr{B}\left(\mathcal{G}\right)\approx\mathcyr{SH}^{1}\left(k,\widehat{\bar{H}}\right)$ (c'est une cons{\'e}quence imm{\'e}diate de l'isomorphisme $\left(3\right)$).
\begin{pro}[Cas global] Soient $k$ un corps de nombres et $\mathcal{G}$ une $k$-gerbe. Pour toute pr{\'e}sentation $\left(V,\pi\right)$ de $\mathcal{G}$:
	\[\mathcyr{B}\left(V\right)\stackrel{\sim}{\longleftarrow}\mathcyr{B}\left(\mathcal{G}\right)
\]
et $m_{\mathcal{H}}\left(\mathcal{G}\right)$ est {\'e}gale {\`a} l'image de $m_{\mathcal{H}}\left(V\right)$ par l'isomorphisme:
	\[\mathcyr{B}\left(V\right)^{D}\stackrel{\sim}{\longrightarrow}\mathcyr{B}\left(\mathcal{G}\right)^{D}
\]
\end{pro}
$\ $
\vspace{1mm}
\newline

L'isomorphisme $\mathcyr{B}\left(V\right)\stackrel{\sim}{\longleftarrow}\mathcyr{B}\left(\mathcal{G}\right)$ r{\'e}sulte de l'isomorphisme 
	\[\mathcyr{SH}^{1}\left(k,Pic\:\bar{V}\right)\stackrel{\sim}{\longrightarrow}\mathcyr{SH}^{1}\left(k,Pic\:\bar{\mathcal{G}}\right)
\]
ce dernier {\'e}tant induit par les isomorphismes compos{\'e}s
	\[H^{1}\left(k,Pic\:\bar{V}\right)\approx H^{1}\left(k,\widehat{\bar{H}}\right)\approx H^{1}\left(k,Pic\:\bar{\mathcal{G}}\right)
\]
Nous avons le diagramme commutatif:
	\[\xymatrix@C=2pt@R=6pt{\prod{V\left(k_{v}\right)}\ar @<-1ex> @{->>}[dd]\ \  \times & \mathcyr{B}\left(V\right) \ar[dr] \\&&& \mathbb{Q}/\mathbb{Z}\\\prod{\mathcal{G}\left(k_{v}\right)}\ \  \times & \mathcyr{B}\left(\mathcal{G}\right) \ar@<1ex>[uu]^{\approx} \ar[ur]}
\]
dans lequel la surjectivit{\'e} de la fl{\`e}che de gauche provient de la proposition 2.1(iv). On a vu que le
calcul de $m_{\mathcal{H}}\left(\mathcal{G}\right)$ ne d{\'e}pendait pas de la famille $\left(x_{v}\right)_{v}$
choisie dans \(\displaystyle\prod_{all\:v}{\mathcal{G}\left(k_{v}\right)}\). De la m{\^e}me mani{\`e}re, on
sait que le calcul de $m_{\mathcal{H}}\left(V\right)$ ne d{\'e}pend pas non plus de la famille
$\left(y_{v}\right)_{v}$ choisie dans \(\displaystyle\prod_{all\:v}{V\left(k_{v}\right)}\). Pour calculer
$m_{\mathcal{H}}\left(V\right)$, on peut donc prendre pour $\left(y_{v}\right)_{v}$ n'importe quel rel{\`e}vement de la famille $\left(x_{v}\right)_{v}$. On en d{\'e}duit que $m_{\mathcal{H}}\left(\mathcal{G}\right)$ n'est autre que l'application compos{\'e}e:
	\[\mathcyr{B}\left(\mathcal{G}\right)\stackrel{\sim}{\longrightarrow}\mathcyr{B}\left(V\right)\xrightarrow{m_{\mathcal{H}}\left(V\right)}\mathbb{Q}/\mathbb{Z}
\]
o{\`u} $m_{\mathcal{H}}\left(V\right)\in \mathcyr{B}\left(V\right)^{D}$ est
donn{\'e}e par: $b\longmapsto \left\langle
\left(y_{v}\right)_{v},b\right\rangle$.

Dans la suite, nous verrons l'{\'e}l{\'e}ment $m_{\mathcal{H}}\left(V\right)$ comme un {\'e}l{\'e}ment de
$\mathcyr{SH}^{1}\left(k,\widehat{\bar{H}}\right)^{D}$.
\end{section}
\begin{section}{$1/2$-th{\'e}or{\`e}me de Tate-Poitou pour les groupes non-ab{\'e}liens}
$\ $

\begin{theo} Soient $k$ un corps de nombres, $H$ un $k$-groupe fini, $\mathcal{L}_{H}$ un $k$-lien localement repr{\'e}sentable par $H$. L'application
	\[\begin{array}{rccc}m_{\mathcal{H}}: & \mathcyr{SH}^{2}\left(k,\mathcal{L}_{H}\right) & \longrightarrow & \mathcyr{SH}^{1}\left(k,\widehat{\bar{H}}\right)^{D}\\ & \left[\mathcal{G}\right] & \longmapsto & m_{\mathcal{H}}\left(\mathcal{G}\right)\end{array}
\]
o{\`u} $\mathcyr{SH}^{2}\left(k,\mathcal{L}_{H}\right)$ d{\'e}signe l'ensemble des classes d'{\'e}quivalence de gerbes localement li{\'e}es par $H$ admettant partout localement une section (i.e qui sont partout localement neutres), se factorise par
	\[\xymatrix{\mathcyr{SH}^{2}\left(k,\mathcal{L}_{H}\right) \ar[r]^{ab} \ar[dr]_{m_{\mathcal{H}}} & \mathcyr{SH}^{2}\left(k,\frac{\bar{H}}{\left[\bar{H},\bar{H}\right]}\right) \ar[d]_{\approx}\\ & \mathcyr{SH}^{1}\left(k,\widehat{\frac{\bar{H}}{\left[\bar{H},\bar{H}\right]}}\right)^{D}}
\]
o{\`u} $ab$ est l'application d'ab{\'e}lianisation naturelle, et l'isomorphisme vertical est fourni par la dualit{\'e} de Tate-Poitou.
\end{theo}
$\ $

\begin{rem}
\item \textup{On peut {\'e}tendre le th{\'e}or{\`e}me 4.1 au cas o{\`u} $H$ est un $k$-groupe lin{\'e}aire,
\textit{i.e} au cas o{\`u} les $k$-gerbes consid{\'e}r{\'e}es ne sont plus de Deligne-Mumford. Ceci peut de
faire en rempla\c{c}ant dans la construction fondamentale (de la section 2) $SL_{n}$ par $GL_{n}$. Le th{\'e}or{\`e}me 4.1 peut alors {\^e}tre compl{\'e}t{\'e} par les deux r{\'e}sultats suivants:}
\begin{enumerate}
\item \textup{Si $H$ est un $k$-tore $T$, $m_{\mathcal{H}}$ prend ses valeurs dans $\mathcyr{SH}^{1}\left(k,X^{\ast}\left(T\right)\right)^{D}$:
	\[m_{\mathcal{H}}:\mathcyr{SH}^{2}\left(k,T\right)\longrightarrow\mathcyr{SH}^{1}\left(k,X^{\ast}\left(T\right)\right)^{D}
\]
et coïncide avec l'isomorphisme donn{\'e} par la dualit{\'e} de Kottwitz pour les tores \cite{K}.}
\item \textup{Si $H$ est un $k$-groupe semi-simple, alors $\mathcyr{SH}^{1}\left(k,\widehat{\bar{H}}\right)=0$ et 
	\[m_{\mathcal{H}}:\mathcyr{SH}^{2}\left(k,\mathcal{L}_{H}\right)\longrightarrow\mathcyr{SH}^{1}\left(k,\widehat{\bar{H}}\right)^{D}=0
\]
est l'application nulle. On sait d'apr{\`e}s \cite{Bo2} dans le cas semi-simple (resp. d'apr{\`e}s \cite{D1} dans le cas semi-simple simplement connexe) que toutes les classes de $\mathcyr{SH}^{2}\left(k,\mathcal{L}_{H}\right)$ (resp. de $H^{2}\left(k,\mathcal{L}_{H}\right)$)
sont neutres. Ainsi l'obstruction de Brauer-Manin est la seule dans le cas
semi-simple. Compte tenu de la remarque $\left(1\right)$ pr{\'e}c{\'e}dente, on
en d{\'e}duit que le m{\^e}me r{\'e}sultat vaut dans le cas des groupes
r{\'e}ductifs connexes, puis dans le cas des groupes connexes (cf \cite{Bo1}).}
\end{enumerate}
\end{rem}
\end{section}

\vspace{20mm}
\[\begin{array}{lll}\textup{\footnotesize{douai@agat.univ-lille1.fr}}\\\textup{\footnotesize{emsalem@agat.univ-lille1.fr}}\\\textup{\footnotesize{stephane.zahnd@agat.univ-lille1.fr}}\end{array}
\]
\end{document}